\documentclass[3p,12pt]{elsarticle}

\makeatletter
\def\ps@pprintTitle{%
\let\@oddhead\@empty
\let\@evenhead\@empty
\def\@oddfoot{}%
\let\@evenfoot\@oddfoot}
\makeatother

\usepackage{amssymb,amsthm,amsmath,natbib}
\usepackage[colorlinks]{hyperref}

\usepackage{geometry}
\newtheorem{theorem}{Theorem}

\DeclareRobustCommand{\brkbinom}{\genfrac\{\}{0pt}{}}
\newcommand{\N}{\mathbb{N}}

\newcommand{\bb}[1]{\boldsymbol{#1}}
\newcommand{\OO}{\mathcal{O}}

\newcommand{\EE}{\mathsf{E}} 

\begin{document}

\begin{frontmatter}

\title{The Gaussian product inequality conjecture\\for multinomial covariances}

\author[a1]{Fr\'ed\'eric Ouimet}\ead{frederic.ouimet@umontreal.ca}

\address[a1]{Centre de recherches math\'ematiques, Universit\'e de Montr\'eal, Montr\'eal (Qu\'ebec) Canada H3T 1J4\vspace{-5mm}}

\begin{abstract}
In this short note, we find an equivalent combinatorial condition only involving finite sums under which a centered Gaussian random vector with multinomial covariance matrix satisfies the Gaussian product inequality (GPI) conjecture. These covariance matrices are relevant since their off-diagonal elements are negative, which is the hardest case to cover for the GPI conjecture, as mentioned by [Russell, O., \& Sun, W. 2022. Some new {G}aussian product inequalities. {\em J. Math. Anal. Appl.}, {\bf 515}(2), Paper No. 126439, 21 pp.].
\end{abstract}

\begin{keyword}
Gaussian product inequality \sep Gaussian random vector \sep moments \sep multinomial distribution \sep multivariate normal \sep stochastic order
\MSC[2020]{Primary: 60E15; Secondary: 05A20, 62E20, 62H10}
\end{keyword}

\end{frontmatter}

\vspace{-4mm}
\section{Introduction}\label{sec:intro}

The Gaussian product inequality (GPI) is a long-standing conjecture attributed to \citet{MR1636556}, which is originally framed in the context of the real polarization problem in functional analysis. The connection to joint moment inequalities for Gaussian random vectors was underscored by \citet{MR2385646}. In this latter framework, the GPI states that for any centered Gaussian random vector $\bb{Z} = (Z_1, \dots, Z_d)$ of dimension $d \in \mathbb{N} = \{ 1, 2, \ldots \}$ and every integer $m \in \mathbb{N}_0 = \{0,1,\ldots\}$, one~has
\begin{equation}\label{eq:1}
\EE \left(\prod_{i=1}^d Z_i^{2 m}\right) \geq \prod_{i=1}^d \EE \big( Z_i^{2 m} \big).
\end{equation}
\citet{MR3425898} proved that inequality~\eqref{eq:1} implies the validity of the real polarization constant conjecture. It is also related to the $U$-conjecture, which asserts that if $P$ and $Q$ are two distinct polynomials on $\mathbb{R}^d$ and the random variables $P(\bb{Z})$ and $Q(\bb{Z})$ are independent, then there exist an orthogonal transformation $L$ on $\mathbb{R}^d$ and an integer $d_1 \in \{ 1, \ldots , d - 1 \}$ such that $P \circ L$ is a function of $(Z_1, \ldots, Z_{d_1})$ and $Q \circ L$ is a function of $(Z_{d_1+1}, \ldots, Z_d)$; see, e.g., \citet{MR0346969,MR3425898}, and references therein.

The inequality \eqref{eq:1} holds for $m = 1$ across every dimension $d \in \N$, a result demonstrated by \citet{MR2385646}. The inequality also holds when $|\bb{Z}| = (|Z_1|, \ldots, |Z_d|)$ possesses the multivariate total positivity of order $2$ (denoted $\mathrm{MTP}_2$) over the domain $[0,\infty)^d$; refer to Corollary~1.1 in \citet{MR628759} and the associated references. Notably, the validity of inequality~\eqref{eq:1} in dimension $d = 2$ is ensured because a random pair $(Z_1, Z_2)$ being Gaussian implies that $(|Z_1|, |Z_2|)$ is $\mathrm{MTP}_2$, as mentioned in Remark~1.4 of \citet{MR628759}. Building on this, \citet{MR4471184}, employing a moment formula from \citet{MR0045347}, established that the converse of inequality \eqref{eq:1} holds for $d = 2$ with $(\alpha_1,\alpha_2)\in (-1,0] \times [0,\infty)$, thereby completing the exploration of the GPI conjecture for the bivariate case. A more precise version of the GPI, offering tighter bounds, was investigated by \citet{MR4530374}.

Interest in the conjecture was renewed recently by the proof of \citet{MR4052574} that inequality~\eqref{eq:1} holds in dimension $d = 3$, based on a careful analysis of Gaussian hypergeometric functions.

As originally proposed by \citet{MR2886380}, a certain stronger version of the conjecture is believed to be true. The strong version states that for any centered Gaussian random vector $\bb{Z} = (Z_1, \dots, Z_d)$ of dimension $d \in \mathbb{N}$ and every reals $\alpha_1, \ldots, \alpha_d \in [0, \infty)$, one~has
\begin{equation}\label{eq:2}
\EE \left( \prod_{i=1}^d |Z_i|^{2 \alpha_i} \right) \geq \prod_{i=1}^d \EE \big (|Z_i|^{2 \alpha_i} \big).
\end{equation}
For nonnegative integer exponents, \citet{MR4466643} proved inequality~\eqref{eq:2} under the assumption that the covariance matrix $\Sigma$ is completely positive (meaning that there exists a matrix $C$ with nonnegative entries such that $\Sigma = C C^{\top}$). Their proof utilized a combinatorial methodology tied closely to the complete monotonicity of multinomial probabilities, a result previously obtained independently by \citet{MR3825458} and \citet{MR4201158}. The findings of Genest and Ouimet were soon after broadened to encompass all covariance matrices having nonnegative entries, a refinement made by \citet{MR4445681} using the Isserlis--Wick type formula of \citet{MR3324071,arXiv:1705.00163}; further details can be found in Corollary~1 of \citet{arXiv:2202.00189} and Remark~2.3 of \citet{MR4554766}. Similarly, the work of Russell and Sun was extended by \citet{MR4554766} to the case where the variables $Z_1^2, \ldots, Z_d^2$ are replaced by components of a random vector $\bb{X} = (X_1, \ldots, X_d)$ which possesses a multivariate gamma distribution in the sense of \citet{MR44790}. Here, the covariance matrix $\Sigma$ is such that a signature matrix $S$ can be found where every entry of $S\Sigma S$ is nonnegative. \citet{MR4538422} subsequently highlighted that this latter result can be extended even further to account for random vectors of traces of diagonal blocks of Wishart matrices, in line with the authors' belief that a specific Kronecker product version of the GPI should hold more generally for these matrices. The concept of random vectors of traces of diagonal blocks of Wishart matrices was initially discussed by \citet{MR263201} in relation to the renowned Gaussian correlation inequality (GCI) conjecture, which saw its resolution in 2014 by \citet{MR3289621}, see also \citet{MR3468024,MR3645127}.
For a comprehensive review of the properties of the multivariate gamma distribution in the sense of \citet{MR44790}, the interested readers are referred to \citet{MR3325368}.

For the case where $d = 3$ and the exponents are $(\alpha_1, \alpha_2, \alpha_3) \in \{1\} \times \{2, 3\} \times \N_0$, inequality~\eqref{eq:2} was validated by \citet{MR4445681}, who employed a brute-force combinatorial approach. Taking advantage of a sums-of-squares methodology combined with in-depth symbolic/numerical calculations using \texttt{Macaulay2} and \texttt{Mathematica}, \citet{arXiv:2205.02127} asserted in their Theorems~4.1~and~4.2 the validity of inequality~\eqref{eq:2} for the sets $(d, \alpha_1, \alpha_2, \alpha_3) \in \{3\} \times \N_0 \times \{3\} \times \{2\}$ and $(d, \alpha_1, \alpha_2, \alpha_3, \alpha_4) \in \{4\} \times \N_0 \times \{2\} \times \{2\} \times \{2\}$ respectively. A moment ratio inequality related to the three-dimensional GPI was also explored by \citet{MR4593134}.

The most important recent development is due to \citet{arXiv:2211.07314}, who showed, using induction and an elegant Gaussian integration by parts argument, that inequality~\eqref{eq:2} holds in dimension $d = 3$ for all nonnegative integer exponents. It remains to be seen if their approach can be extended to cover all dimensions.

Variants of the weak and strong versions of the GPI are also considered in the literature. For instance, \citet{MR3278931} proved that, for the aforementioned multivariate gamma distributed random vector $\bb{X} = (X_1,\ldots,X_d)$ in the sense of \citet{MR44790}, together with a certain range of negative exponents $\beta_1, \ldots, \beta_d$ which depends on the degree of freedom, and any integer $d_1\in \{1, \ldots, d-1 \}$, one has
\begin{equation}\label{eq:3}
\EE \left(\prod_{i=1}^d X_i^{\beta_i} \right) \geq \EE \left(\prod_{i=1}^{d_1} X_i^{\beta_i}\right) \EE \left(\prod_{i=d_1+1}^d X_i^{\beta_i}\right).
\end{equation}
The special case where the $X_i$'s are the squared components of a Gaussian random vector holds for any negative exponents $\beta_i = 2\alpha_i\in (-1,0), ~i\in \{1,\ldots,d\}$. Inequality~\eqref{eq:3} was recovered by \citet{MR4554766} through the application of the multivariate gamma extension of the GCI, due to \citet{MR3289621}. These researchers noted that the component-wise absolute negative powers of multivariate gamma random vectors are strongly positive upper orthant dependent, which then allowed them to integrate on both sides of the associated inequality to obtain \eqref{eq:3} for negative exponents (whenever the moments are finite).

Through the application of the multivariate Laplace transform order (see, e.g., \citet[Section~7.D.1]{MR2265633}), the even more general inequality
\begin{equation}\label{eq:4}
\EE\left\{\prod_{i=1}^d \phi_i(X_i)\right\} \geq \EE\left\{\prod_{i=1}^{d_1} \phi_i(X_i)\right\} \EE\left\{\prod_{i=d_1+1}^{d} \phi_i(X_i)\right\},
\end{equation}
was shown by \citet{MR4538422}, where $\bb{X} = (X_1,\ldots,X_d)$ is a random vector of traces of diagonal blocks inside a Wishart matrix, and the functions $\phi_i$'s are only assumed to be completely monotonic on $(0,\infty)$ (note that the functions $\phi_i(x_i) = x_i^{\beta_i}$ are completely monotonic on $(0,\infty)$ for any $\beta_i < 0$). In particular, Corollary~1 of that paper showed how \eqref{eq:4} generalizes multiple GPI-type inequalities derived earlier by \citet{MR3278931} and \citet{MR3608204}. It also introduced an intriguing relationship with the famous BMV inequality (formerly known as the BMV conjecture) proved by \citet{MR3143891}, which states that for any symmetric positive definite matrix $A$ of size $n\times n$, the map $t \mapsto \mathrm{tr} \{ \exp(-t A) \}$ is completely monotonic on $(0,\infty)$.

In this short note, our objective is to delve into the GPI inequality~\eqref{eq:1}, focusing on a multivariate Gaussian random vector $\bb{Z}$ whose covariance matrix adopts a multinomial structure. In Theorem~\ref{thm:1}, we identify a combinatorial condition equivalent to the original inequality, characterized solely by finite sums, offering a numerically attractive perspective. As pointed out by \citet{MR4445681}, such covariance matrices are of particular interest; their off-diagonal elements are negative, representing the most challenging scenario when studying the GPI conjecture.

Section~\ref{sec:multinomial} introduces the multinomial distribution and its covariance matrix. The aforementioned GPI equivalence is stated in Section~\ref{sec:result} and proved in Section~\ref{sec:proof}.

\section{The multinomial distribution}\label{sec:multinomial}

For any $d\in \N$, the interior of the $d$-dimensional (unit) simplex is defined by
\[
\mathcal{S}_d^o = \big\{\textstyle\bb{x}\in (0,1)^d : \sum_{i=1}^d x_i\in (0,1)\big\}.
\]
A random vector $\bb{\xi} = (\xi_1,\xi_2,\dots,\xi_d) \sim \mathrm{Multinomial}\hspace{0.2mm}(N,\bb{p})$ follows a multinomial distribution if its probability mass function is defined by
\begin{equation}\label{eq:5}
p_{\bb{\xi}}(\bb{k}) = \frac{N!}{(N - \sum_{i=1}^d k_i)! \prod_{i=1}^d k_i!} (1 - \sum_{i=1}^d p_i)^{N - \sum_{i=1}^d k_i} \prod_{i=1}^d p_i^{k_i}, \quad \bb{k}\in \N_0^d \cap N \mathcal{S}_d^o,
\end{equation}
where $N\in \N$ and $\bb{p}\in \mathcal{S}_d^o$ are parameters. The covariance matrix of the multinomial distribution is well known to be
\vspace{-1mm}
\[
\Sigma = \mathrm{diag}(\bb{p}) - \bb{p} \bb{p}^{\top},
\]
see, e.g., \citet{MR1157720} or Chapter~35 of \citet{MR1429617}.

\section{Result}\label{sec:result}

Let $\bb{\xi} \sim \mathrm{Multinomial}\hspace{0.2mm}(N,\bb{p})$ and $\bb{Z} \sim \mathcal{N}_d(\bb{0}_d, \Sigma)$. By the central limit theorem,
\begin{equation}\label{eq:6}
\frac{\bb{\xi} - N \bb{p}}{\sqrt{N}} = \left(\frac{\xi_i - N p_i}{\sqrt{N}}\right)_{i=1}^d \stackrel{\mathrm{law}}{\longrightarrow} \bb{Z}, \quad \text{as } N\to \infty,
\end{equation}
see, e.g., \cite{MR0478288,MR750392,MR1157720,MR4249129,MR4361955} for closely related local approximations. Consequently, the GPI~\eqref{eq:1} holds for multinomial covariances if and only if
\begin{equation}\label{eq:7}
\lim_{N\to \infty} \left[\EE\left\{\prod_{i=1}^d \left(\frac{\xi_i - N p_i}{\sqrt{N}}\right)^{2m}\right\} - \prod_{i=1}^d \EE\left\{\left(\frac{\xi_i - N p_i}{\sqrt{N}}\right)^{2m}\right\}\right] \geq 0.
\end{equation}
In Theorem~\ref{thm:1} below, we find a combinatorial condition equivalent to~\eqref{eq:7} only involving finite sums.

\begin{theorem}\label{thm:1}
Let $d,m,N\in \N$ and $\bb{p}\in \mathcal{S}_d^o$ be given, and let $\bb{\xi}\sim \mathrm{Multinomial}\hspace{0.2mm}(N,\bb{p})$ be defined as in~\eqref{eq:5}. The inequality
\[
\lim_{N\to \infty} \left[\EE\left\{\prod_{i=1}^d \left(\frac{\xi_i - N p_i}{\sqrt{N}}\right)^{2m}\right\} - \prod_{i=1}^d \EE\left\{\left(\frac{\xi_i - N p_i}{\sqrt{N}}\right)^{2m}\right\}\right] \geq 0
\]
is equivalent to
\begin{equation}\label{eq:8}
\begin{aligned}
&\sum_{\ell_1, \dots, \ell_d=0}^{2m} \sum_{k_1=0}^{\ell_1} \dots \sum_{k_d=0}^{\ell_d} \left\{\left[\hspace{-1mm}\begin{array}{c}\sum_{i=1}^d k_i\\\sum_{i=1}^d \ell_i - md\end{array}\hspace{-1mm}\right]^{\star} - \prod_{i=1}^d \left[\hspace{-1mm}\begin{array}{c}k_i\\\ell_i - m\end{array}\hspace{-1mm}\right]^{\star}\right\} \\
&\quad\hspace{40mm}\times \prod_{i=1}^d \left[\brkbinom{\ell_i}{k_i} \binom{2m}{\ell_i} (-1)^{\ell_i} p_i^{2m - \ell_i + k_i}\right] \geq 0,
\end{aligned}
\end{equation}
where {\scriptsize $\Big[\hspace{-1.3mm}\begin{array}{c}k\\[-0.8mm]j\end{array}\hspace{-1.3mm}\Big]^{\hspace{-0.5mm}\star}$} \hspace{-1mm}denotes a signed Stirling number of the first kind and $\brkbinom{k}{j}$ denotes a Stirling number of the second kind, with the understanding that these two quantities are zero whenever $j < 0$ or $j > k$.
\end{theorem}

\section{Proof}\label{sec:proof}

The proof of the equivalent condition \eqref{eq:8} uses the falling factorial decomposition of monomials combined with the general expressions for the falling factorial moments of the multinomial distribution obtained by \citet{MR143299}. For $d = m = 2$, the multinomial GPI can be verified directly (for example) by using the central multinomial moment formulas derived in Section~5.2 of \citet{doi:10.3390/stats4010002} or p.\hspace{0.3mm}4--5 of \citet{arXiv:2006.09059}.

\begin{proof}[Proof of Theorem~\ref{thm:1}]
Denote
\[
x^{(0)} = 1, \qquad x^{(j)} = x (x - 1) \dots (x - j + 1), \quad j\in \N.
\]
From Equations~(6.10)~and~(6.13) of \citet{MR1397498}, we know that falling factorials can be expressed as specific linear combinations of monomials and conversely:
\begin{equation}\label{eq:9}
x^{(k)} = \sum_{j=0}^k \left[\hspace{-1mm}\begin{array}{c}k\\j\end{array}\hspace{-1mm}\right]^{\star} x^j, \qquad x^{\ell} = \sum_{k=0}^{\ell} \brkbinom{\ell}{k} x^{(k)}, \qquad k,\ell\in \N_0.
\end{equation}
(Compared to (6.13) of \citet{MR1397498}, we have ${\scriptsize \Big[\hspace{-1.3mm}\begin{array}{c}k\\[-0.8mm]j\end{array}\hspace{-1.3mm}\Big]^{\star}} = (-1)^{k-j} {\scriptsize \Big[\hspace{-1.3mm}\begin{array}{c}k\\[-0.8mm]j\end{array}\hspace{-1.3mm}\Big]}$, where ${\scriptsize \Big[\hspace{-1.3mm}\begin{array}{c}k\\[-0.8mm]j\end{array}\hspace{-1.3mm}\Big]}$ denotes a unsigned Stirling number of the first kind.) Also, from \citet{MR143299}, we have explicit expressions for the falling factorial moments of the multinomial distribution:
\begin{equation}\label{eq:10}
\EE\left\{\xi_1^{(k_1)} \dots \xi_d^{(k_d)}\right\} = N^{(\sum_{i=1}^d k_i)} \prod_{i=1}^d p_i^{k_i} \quad \text{and} \quad \EE\left\{\xi_i^{(k_i)}\right\} = N^{(k_i)} p_i^{k_i}.
\end{equation}

By applying the binomial formula together with~\eqref{eq:9} and~\eqref{eq:10}, we have
\[
\begin{aligned}
&\EE\left\{\prod_{i=1}^d \left(\frac{\xi_i - N p_i}{\sqrt{N}}\right)^{2m}\right\} \notag \\
&= N^{-md} \sum_{\ell_1, \dots, \ell_d=0}^{2m} \EE\left\{\xi_1^{\ell_1} \dots \xi_d^{\ell_d}\right\} \prod_{i=1}^d \left[\binom{2m}{\ell_i} (-N p_i)^{2m - \ell_i}\right] \notag \\
&= N^{-md} \sum_{\ell_1, \dots, \ell_d=0}^{2m} \sum_{k_1=0}^{\ell_1} \dots \sum_{k_d=0}^{\ell_d} \EE\left\{\xi_1^{(k_1)} \dots \xi_d^{(k_d)}\right\} \prod_{i=1}^d \left[\brkbinom{\ell_i}{k_i} \binom{2m}{\ell_i} (-N p_i)^{2m - \ell_i}\right] \notag \\
&= N^{-md} \sum_{\ell_1, \dots, \ell_d=0}^{2m} \sum_{k_1=0}^{\ell_1} \dots \sum_{k_d=0}^{\ell_d} N^{(\sum_{i=1}^d k_i)} \prod_{i=1}^d \left[\brkbinom{\ell_i}{k_i} \binom{2m}{\ell_i} (-N)^{2m - \ell_i} p_i^{2m - \ell_i + k_i}\right] \notag \\
&= \sum_{\ell_1, \dots, \ell_d=0}^{2m} \sum_{k_1=0}^{\ell_1} \dots \sum_{k_d=0}^{\ell_d} \sum_{j=0}^{\sum_{i=1}^d k_i} N^{md - \sum_{i=1}^d \ell_i + j} \left[\hspace{-1mm}\begin{array}{c}\sum_{i=1}^d k_i\\j\end{array}\hspace{-1mm}\right]^{\star} \notag \\[-1mm]
&\hspace{55mm}\times \prod_{i=1}^d \left[\brkbinom{\ell_i}{k_i} \binom{2m}{\ell_i} (-1)^{\ell_i} p_i^{2m - \ell_i + k_i}\right].
\end{aligned}
\]
Note that the contributions coming from the terms satisfying
\[
md - \sum_{i=1}^d \ell_i + j \neq 0
\]
are forced to be $0$ or $\OO(N^{-1})$ because of the asymptotic normality in~\eqref{eq:6} and the fact that all Gaussian moments of the form $\EE(\prod_{i=1}^d Z_i^{2m})$ are finite and non-zero. Therefore, only the term $j = \sum_{i=1}^d \ell_i - md$ matters in the limit, so that
\begin{equation}\label{eq:11}
\begin{aligned}
&\lim_{N\to \infty} \EE\left\{\prod_{i=1}^d \left(\frac{\xi_i - N p_i}{\sqrt{N}}\right)^{2m}\right\} \\
&= \sum_{\ell_1, \dots, \ell_d=0}^{2m} \sum_{k_1=0}^{\ell_1} \dots \sum_{k_d=0}^{\ell_d} \left[\hspace{-1mm}\begin{array}{c}\sum_{i=1}^d k_i\\\sum_{i=1}^d \ell_i - md\end{array}\hspace{-1mm}\right]^{\star} \prod_{i=1}^d \left[\brkbinom{\ell_i}{k_i} \binom{2m}{\ell_i} (-1)^{\ell_i} p_i^{2m - \ell_i + k_i}\right].
\end{aligned}
\end{equation}

Similarly, by applying the binomial formula together with~\eqref{eq:9} and~\eqref{eq:10}, we have, for all $i\in \{1,2,\dots,d\}$,
\[
\begin{aligned}
\EE\left\{\left(\frac{\xi_i - N p_i}{\sqrt{N}}\right)^{2m}\right\}
&= N^{-m} \sum_{\ell_i=0}^{2m} \EE\left\{\xi_i^{\ell_i}\right\} \binom{2m}{\ell_i} (-N p_i)^{2m - \ell_i} \notag \\
&= N^{-m} \sum_{\ell_i=0}^{2m} \sum_{k_i=0}^{\ell_i} \EE\left\{\xi_i^{(k_i)}\right\} \brkbinom{\ell_i}{k_i} \binom{2m}{\ell_i} (-N p_i)^{2m - \ell_i} \notag \\
&= N^{-m} \sum_{\ell_i=0}^{2m} \sum_{k_i=0}^{\ell_i} N^{(k_i)} \brkbinom{\ell_i}{k_i} \binom{2m}{\ell_i} (-N)^{2m - \ell_i} p_i^{2m - \ell_i + k_i} \notag \\
&= \sum_{\ell_i=0}^{2m} \sum_{k_i=0}^{\ell_i} \sum_{j_i=0}^{k_i} N^{m - \ell_i + j_i} \left[\hspace{-1mm}\begin{array}{c}k_i\\j_i\end{array}\hspace{-1mm}\right]^{\star} \brkbinom{\ell_i}{k_i} \binom{2m}{\ell_i} (-1)^{\ell_i} p_i^{2m - \ell_i + k_i}.
\end{aligned}
\]
Note that the contributions coming from the terms satisfying $m - \ell_i + j_i \neq 0$ are forced to be $0$ or $\OO(N^{-1})$ because the central limit theorem shows that
\[
\frac{\xi_i - N p_i}{\sqrt{N}} \stackrel{\mathrm{law}}{\longrightarrow} Z_i \sim \mathcal{N}(0, p_i - p_i^2), \quad \text{as } N\to \infty,
\]
and all Gaussian moments of the form $\EE(Z_i^{2m})$ are finite and non-zero. Therefore, only the terms $j_i = \ell_i - m$ matter in the limit, so that, for all $i\in \{1,2,\dots,d\}$,
\[
\lim_{N\to \infty} \EE\left\{\left(\frac{\xi_i - N p_i}{\sqrt{N}}\right)^{2m}\right\} = \sum_{\ell_i=0}^{2m} \sum_{k_i=0}^{\ell_i} \left[\hspace{-1mm}\begin{array}{c}k_i\\\ell_i - m\end{array}\hspace{-1mm}\right]^{\star} \brkbinom{\ell_i}{k_i} \binom{2m}{\ell_i} (-1)^{\ell_i} p_i^{2m - \ell_i + k_i}.
\]
If we take the product in the index variable $i$, we get
\begin{equation}\label{eq:12}
\begin{aligned}
&\lim_{N\to \infty} \prod_{i=1}^d \EE\left\{\left(\frac{\xi_i - N p_i}{\sqrt{N}}\right)^{2m}\right\} \\
&= \sum_{\ell_1, \dots, \ell_d=0}^{2m} \sum_{k_1=0}^{\ell_1} \dots \sum_{k_d=0}^{\ell_d} \prod_{i=1}^d \left[\hspace{-1mm}\begin{array}{c}k_i\\\ell_i - m\end{array}\hspace{-1mm}\right]^{\star} \prod_{i=1}^d \left[\brkbinom{\ell_i}{k_i} \binom{2m}{\ell_i} (-1)^{\ell_i} p_i^{2m - \ell_i + k_i}\right].
\end{aligned}
\end{equation}
The conclusion follows by subtracting \eqref{eq:12} from \eqref{eq:11}.
\end{proof}

\section*{Funding}

F.\ Ouimet is supported by a CRM-Simons postdoctoral fellowship from the Centre de recherches math\'ematiques (Montr\'eal, Canada) and the Simons Foundation.

\addcontentsline{toc}{section}{References}

\bibliographystyle{authordate1}
\bibliography{Ouimet_2023_GPI_multinomial_bib}

\end{document}